\documentclass[12pt]{article}
\usepackage{amsmath}
\usepackage[affil-it]{authblk}
\usepackage{amsfonts}
\usepackage{amsthm}
\usepackage{amssymb}
\usepackage{makecell}
\usepackage{appendix}
\usepackage{algorithmic}
\usepackage{graphicx}
\usepackage{algorithm2e}
\usepackage{tabu}
\usepackage{color}
\definecolor{darkgreen}{rgb}{0.0, 0.63, 0.0}
\theoremstyle{definition}

\usepackage{mathtools}

\theoremstyle{lemma}
\newtheorem{theorem}{Theorem}

\newtheorem{proposition}{Proposition}

\usepackage[margin=0.5in]{geometry}

\begin{document}
    \title{A New Bound for the Orthogonality Defect of HKZ Reduced Lattices}
    \author{Christian Porter, Edmund Dable-Heath, Cong Ling}
    \maketitle
\begin{abstract}
    In this work, we determine a sharp upper bound on the orthogonality defect of HKZ reduced bases up to dimension $3$. Using this result, we determine a general upper bound for the orthogonality defect of HKZ reduced bases of arbitrary rank. This upper bound seems to be sharper than existing bounds in literature, such as the one determined by Lagarias, Lenstra and Schnorr \cite{HKZ}.
\end{abstract}
\section{Introduction}
A lattice $\Lambda$ is a discrete subgroup of $\mathbb{R}^m$, for some positive integer $m$. Every lattice $\Lambda$ has a basis $B=\{\mathbf{b}_1,\mathbf{b}_2,\dots, \mathbf{b}_n\}$, and each point of $\Lambda$ may be represented as the linear sum of its basis vectors over $\mathbb{Z}$, that is,
\begin{align*}
    \Lambda=\Lambda(B)=\left\{\sum_{i=1}^n x_i \mathbf{b}_i: x_i \in \mathbb{Z}\right\}.
\end{align*}
We say that $\Lambda$ is of rank $n$ if $\mathbf{b}_1,\dots,\mathbf{b}_n$ are linearly independent over $\mathbb{R}$ and form a basis for $\mathbb{R}^n$, and we say that $\Lambda$ is full-rank if $n=m$.
\\
The process of moving from a ``bad'' lattice basis (one in which the vectors are relatively long, non-orthogonal) to a ``good'' basis (relatively short and orthogonal vectors) is known as \emph{reduction theory}. There are many different definitions of what constitutes a reduced lattice basis, and the focus of this paper will be on so-called HKZ reduced bases, named for mathematicians C. Hermite, A. Korkin and G. Zolotarev \cite{hermitered} \cite{surlesformesquadratiques}. Denote by
\begin{align*}
    &\mathbf{b}_i(j)=\mathbf{b}_i-\sum_{k=1}^{j-1}\mu_{i,k}\mathbf{b}_k(k), \hspace{2mm} \forall 1\leq j \leq i \leq n,
    \\&\mathbf{b}_1(1)=\mathbf{b}_1,
    \\& \mu_{i,j}=\frac{\langle \mathbf{b}_i,\mathbf{b}_j(j)\rangle}{\|\mathbf{b}_j(j)\|^2}, \hspace{2mm} \forall 1 \leq j <i \leq n.
\end{align*}
Then $B=\{\mathbf{b}_1,\mathbf{b}_2,\dots,\mathbf{b}_n\}$ is said to be HKZ reduced if the following properties hold:
\begin{itemize}
    \item $|\mu_{i,j}| \leq 1/2$, for all $1 \leq j < i \leq n$,
    \item $\mathbf{b}_1$ is the shortest nonzero lattice vector of $\Lambda(B)$,
    \item $\{\mathbf{b}_2(2),\mathbf{b}_3(2),\dots,\mathbf{b}_n(2)\}$ is HKZ reduced.
\end{itemize}
The \emph{$i$th successive minimum} of a lattice $\Lambda$, denoted by $\lambda_i$, is the smallest real number such that there are $i$ linearly independent vectors in $\Lambda$ of length at most $\lambda_i$.
\begin{proposition}[\cite{surlesformesquadratiques}] \label{propn1}
    If $B=\{\mathbf{b}_1,\mathbf{b}_2,\dots,\mathbf{b}_n\}$ is HKZ reduced, then for all $1 \leq i \leq n-1$,
    \begin{align*}
        \|\mathbf{b}_i(i)\|^2 \leq \frac{4}{3}\|\mathbf{b}_{i+1}(i+1)\|^2,
    \end{align*}
    and for all $1 \leq i \leq n-2$,
    \begin{align*}
        \|\mathbf{b}_i(i)\|^2 \leq \frac{3}{2}\|\mathbf{b}_{i+2}(i+2)\|^2.
    \end{align*}
\end{proposition}
\begin{proposition}[\cite{HKZ}]
    If $B=\{\mathbf{b}_1,\mathbf{b}_2,\dots,\mathbf{b}_n\}$ is HKZ reduced, then for all $1 \leq i \leq n$,
    \begin{align*}
        \frac{4}{i+3} \lambda_i^2 \leq \|\mathbf{b}_i\|^2 \leq \frac{i+3}{4} \lambda_i^2.
    \end{align*}
\end{proposition}
We define the Hermite invariant of a lattice $\gamma(B)$ by
\begin{align*}
    \gamma(B)=\frac{\lambda_1^2}{|\det(B)|^{\frac{2}{n}}},
\end{align*}
and the Hermite constant of rank $n$, denoted by $\gamma_n$, is the supremum of all Hermite invariants for rank $n$ lattice bases.

As mentioned, the process of reduction is attaining a basis for a lattice that has ``desirable'' properties. One way of measuring the quality of a basis is by evaluating its \emph{orthogonality defect}, which is defined by
\begin{align*}
    \Delta(B)=\frac{\prod_{i=1}^n \|\mathbf{b}_i\|^2}{\det(B)^2}.
\end{align*}
It is easily seen that the orthogonality defect is not unique to a lattice, and a basis can be transformed so that the orthogonality defect becomes arbitrarily large. However, it is known that the orthogonality defect can be bounded if the basis is HKZ reduced.
\begin{theorem}[\cite{HKZ}]\label{thm2}
    If $B=\{\mathbf{b}_1,\mathbf{b}_2,\dots,\mathbf{b}_n\}$ is HKZ reduced, then
    \begin{align*}
        \Delta(B) \leq \gamma_n^n \prod_{i=1}^n \frac{i+3}{4}.
    \end{align*}
\end{theorem}
Whilst Lagarias et. al. managed to show that the orthogonality defect of HKZ reduced bases must be bounded, little is known about the exact value of the supremum of the orthogonality defect for HKZ reduced basis of a given rank. In this work, we calculate this supremum up to dimension 3 HKZ reduced bases, and use this result to improve the bound in Theorem \ref{thm2}.
\begin{theorem}\label{thm3}
    Let $H_n$ denote the space of rank $n$ HKZ reduced lattice bases. Set
    \begin{align*}
            \Delta_n=\max_{B \in H_n} \Delta(B).
    \end{align*}
    Then $\Delta_1=1$, $\Delta_2=\frac{4}{3}$, $\Delta_3=\frac{25}{12}$.
\end{theorem}
\begin{theorem}\label{thm4}
    Using the notation as before, for all $n \geq 4$,
    \begin{align*}
        \Delta_n \leq \frac{25}{12} \gamma_{n-3}^{n-3} \prod_{i=4}^n \left(\frac{i}{4}+\frac{29}{24}\right).
    \end{align*}
\end{theorem}
\section{Proof of Theorem \ref{thm3}}
The value of $\Delta_1$ is trivial. For $\Delta_2$, for any basis $B=\{\mathbf{b}_1,\mathbf{b}_2\}$,
\begin{align*}
    \Delta(B)=\frac{\|\mathbf{b}_2\|^2}{\|\mathbf{b}_2(2)\|^2}=\frac{1}{1-\frac{\langle \mathbf{b}_2,\mathbf{b}_1\rangle^2}{\|\mathbf{b}_1\|^2\|\mathbf{b}_2\|^2}} \leq \frac{1}{1-\mu_{2,1}^2} \leq \frac{4}{3},
\end{align*}
and this upper bound is attained by the $A_2$ lattice.
\\
We now draw our attention to the case where $n=3$. We will assume without loss of generality that $\|\mathbf{b}_1\|^2=1$ (the value of $\|\mathbf{b}_1\|^2$ will not affect the orthogonality defect of the lattice), and use the following notations:
    \begin{align*}
        \|\mathbf{b}_2(2)\|^2=k, \|\mathbf{b}_3(3)\|^2=l, \mu_{21}=\lambda, \mu_{31}=\mu, \mu_{32}=\sigma.
    \end{align*}
    Then,
    \begin{align*}
        \|x\mathbf{b}_1+y\mathbf{b}_2+z\mathbf{b}_3\|^2=(x+\lambda y +\mu z)^2+k(y+\sigma z)^2 +lz^2.
    \end{align*}
    By switching the signs of $y,z$ we may assume without loss of generality that $\lambda,\mu \geq 0$. It is easily seen that we have,
    \begin{align*}
        &\frac{\|\mathbf{b}_2\|^2}{\|\mathbf{b}_2(2)\|^2}=1+\frac{\lambda^2}{k},
        \\&\frac{\|\mathbf{b}_3\|^2}{\|\mathbf{b}_3(3)\|^2}=1+\frac{1}{l}(\mu^2+k\sigma^2).
    \end{align*}
    Since the basis is HKZ reduced, the following inequalities hold:
    \begin{align}
        &k+\lambda^2 \geq 1, \label{1}
        \\&l+k\sigma^2+\mu^2 \geq 1, \label{2}
        \\&
        l+k(1+\sigma)^2+(1-\lambda-\mu)^2 \geq 1, \label{3}
        \\&
        l+k(1-\sigma)^2+(\lambda-\mu)^2 \geq 1 \label{4}
        \\&
        l+k\sigma^2 \geq k. \label{5}
    \end{align}
    From this, we may ascertain the following inequalities:
    \begin{align}
        &k \leq \frac{l}{1-\sigma^2}, \label{6}
        \\& \lambda^2 \geq 1-\frac{l}{1-\sigma^2}, \label{7}
        \\& \mu^2 \geq 1-\frac{l}{1-\sigma^2}. \label{8}
    \end{align}
    Now, using the inequalities $l \geq \frac{2}{3},k \geq \frac{3}{4}$, if $|\sigma| \leq 1/3$ we have
    \begin{align*}
        \Delta(B)=\left(1+\frac{\lambda^2}{k}\right)\left(1+\frac{1}{l}(\mu^2+k\sigma^2)\right) \leq\left(1+\frac{1}{4k}\right)\left(1+\frac{\frac{1}{4}+\frac{9}{8}\frac{1}{9}l}{l}\right)=\left(1+\frac{1}{4k}\right)\left(\frac{9}{8}+\frac{1}{4l}\right)\leq \frac{4}{3}\frac{3}{2}=\gamma_3^3,
    \end{align*}
    so we may assume without loss of generality that $1/3 \leq |\sigma| \leq 1/2$. Assume first that $\sigma \leq -1/3$. By inequalities \ref{3} and \ref{6}, we have
    \begin{align*}
        l \geq \max\{1-(1-\lambda-\mu)^2-k(1+\sigma)^2,k(1-\sigma^2)\}.
    \end{align*}
    Assume first that
    \begin{align}
        1-(1-\lambda-\mu)^2-k(1+\sigma)^2 \geq k(1-\sigma^2) \iff k \leq \frac{1}{2(1+\sigma)}\left(1-(1-\lambda-\mu)^2\right). \label{9}
    \end{align}
    Then
    \begin{align*}
        \Delta(B) \leq \left(1+\frac{\lambda^2}{k}\right)\left(1+\frac{\mu^2+k\sigma^2}{1-(1-\lambda-\mu)^2-k(1+\sigma)^2}\right).
    \end{align*}
    Let
    \begin{align*}
        f(k)=\left(1+\frac{\lambda^2}{k}\right)\left(1+\frac{\mu^2+k\sigma^2}{1-(1-\lambda-\mu)^2-k(1+\sigma)^2}\right),
    \end{align*}
    so
    \begin{align*}
        f^{\prime \prime}(k)=\frac{\alpha}{k^3(1-(1-\lambda-\mu)^2-k(1+\sigma)^2)^3},
    \end{align*}
    where
    \begin{align*}
        \alpha&=2\lambda^2(\mu^2+k\sigma^2)(1-(1-\lambda-\mu)^2-k(1+\sigma)^2)^2+2\lambda^2(1-(1-\lambda-\mu)^2-k(1+\sigma)^2)^3\\&-2\lambda^2k(1+\sigma)^2(\mu^2+k\sigma^2)(1-(1-\lambda-\mu)^2-k(1+\sigma)^2)-2k\lambda^2\sigma^2(1-(1-\lambda-\mu)^2-k(1+\sigma)^2)^2\\&+2\lambda^2k^2(1+\sigma)^4(\mu^2+k\sigma^2)+2k^3(1+\sigma)^4(\mu^2+k\sigma^2)+2\lambda^2k^2\sigma^2(1+\sigma)^2(1-(1-\lambda-\mu)^2-k(1+\sigma)^2)\\&+2k^3\sigma^2(1+\sigma)^2
        \\&=\lambda^2(\mu^2+k\sigma^2)(1-(1-\lambda-\mu)^2-2k(1+\sigma)^2)^2\\&+\lambda^2(1-(1-\lambda-\mu)^2-k(1+\sigma)^2)(1-(1-\lambda-\mu)^2-k((1+\sigma)^2+\sigma^2))^2
        \\&+\lambda^2(\mu^2+k\sigma^2)(1-(1-\lambda-\mu)^2-k(1+\sigma)^2)^2+\lambda^2(1-(1-\lambda-\mu)^2-k(1+\sigma)^2)^3
        \\&+\lambda^2k^2(1+\sigma)^4(\mu^2+k\sigma^2)+2k^3(1+\sigma)^4(\mu^2+k\sigma^2)\\&+\lambda^2k^2(2\sigma^2(1+\sigma)^2-\sigma^4)(1-(1-\lambda-\mu)^2-k(1+\sigma)^2)+2k^3\sigma^2(1+\sigma)^2 \geq 0,
    \end{align*}
    since
    \begin{align*}
        2\sigma^2(1+\sigma)^2-\sigma^4=\sigma^2(\sigma^2+4\sigma+2)>0
    \end{align*}
    for all $|\sigma|\leq 1/2$, and so $f^{\prime \prime}(k)$ is non-negative for all values of $\lambda,k,\sigma, \mu$ over their respective regions, so the maximum value of $f(k)$ occurs either when $k$ is maximum or when it is minimum. Suppose first that $f(k)$ is maximum when $k$ is minimum. Then
    \begin{align*}
        \Delta(B) &\leq \left(1+\frac{\lambda^2}{1-\lambda^2}\right)\left(1+\frac{\mu^2+\sigma^2(1-\lambda^2)}{1-(1-\lambda-\mu)^2-(1-\lambda^2)(1+\sigma)^2}\right).
    \end{align*}
    Let's assume a contradiction, so
    \begin{align}
        \left(1+\frac{\lambda^2}{1-\lambda^2}\right)\left(1+\frac{\mu^2+\sigma^2(1-\lambda^2)}{1-(1-\lambda-\mu)^2-(1-\lambda^2)(1+\sigma)^2}\right) \geq \frac{25}{12}, \label{10}
    \end{align}
    for some values of $\mu,\sigma,\lambda$. Rearranging gives us the quadratic in $\sigma$:
    \begin{align*}
        &\frac{25}{12}(1-\lambda^2)^2\sigma^2+\left(-2(1-\lambda^2)+\frac{25}{6}(1-\lambda^2)^2\right)\sigma+1-(1-\lambda-\mu)^2-\frac{37}{12}(1-\lambda^2)+\mu^2\\&+\frac{25}{12}(1-\lambda^2)(1-\lambda-\mu)^2+\frac{25}{12}(1-\lambda^2)^2 \geq 0
    \end{align*}
    Denote by $r^+,r^-$ respectively the greater and lesser roots of the equation above. Then we need either
    \begin{align*}
        \sigma \geq r^+, \hspace{2mm} \sigma \leq r^-.
    \end{align*}
    However, we can show that the solutions to this equation yield contradictions. Using \ref{9}, the function is only valid under the boundaries
    \begin{align}
        &0 \leq \lambda \leq 1/2, \label{11}
        \\& 1+\lambda-\sqrt{\lambda^2+2\lambda} \leq \mu \leq 1/2, \label{12}
    \end{align}
    and the function is undefined if $1+\lambda-\sqrt{\lambda^2+2\lambda}>1/2$. Using these boundaries, we can numerically show that there cannot exist a valid solution as we would require $\sigma >-1/3$ for $\sigma \geq r^+$, or $\sigma<-1/2$ for $\sigma \leq r^-$. Hence, the inequality \ref{10} cannot exist under the assumptions made.
    \\
    Now assume that the maximum of $f(k)$ occurs when $k$ is maximum, so
    \begin{align}
        \Delta(B) \leq \left(1+\frac{2\lambda^2(1+\sigma)}{1-(1-\lambda-\mu)^2}\right)\left(1+\frac{\mu^2+\frac{\sigma^2}{2(1+\sigma)}(1-(1-\lambda-\mu)^2)}{\frac{1}{2}(1-\sigma)(1-(1-\lambda-\mu)^2)}\right). \label{13}
    \end{align}
    Once again, we assume that we have a contradiction so
    \begin{align*}
        \left(1+\frac{2\lambda^2(1+\sigma)}{1-(1-\lambda-\mu)^2}\right)\left(1+\frac{\mu^2+\frac{\sigma^2}{2(1+\sigma)}(1-(1-\lambda-\mu)^2)}{\frac{1}{2}(1-\sigma)(1-(1-\lambda-\mu)^2)}\right) \geq \frac{25}{12},
    \end{align*}
    which yields the quadratic
    \begin{align*}
        &\frac{25\lambda^4+100\lambda^3\mu+198\lambda^2\mu^2+100\lambda\mu^3+25\mu^4-100\lambda^3-300\lambda^2\mu-300\lambda\mu^2-100\mu^3+100\lambda^2+200\lambda\mu+100\mu^2}{12}\sigma^2
        \\&-2\left(\lambda^4+2\lambda^3\mu-2\lambda^2\mu^2+2\lambda\mu^3+\mu^4-2\lambda^3-2\lambda^2\mu-2\lambda\mu^2-2\mu^3\right)\sigma
        \\&-\frac{37}{12}\lambda^4-\frac{25}{3}\lambda^3\mu-\frac{13}{2}\lambda^2\mu^2-\frac{25}{3}\lambda\mu^3-\frac{37}{12}\mu^4+\frac{25}{3}\lambda^3+17\lambda^2\mu+17\lambda\mu^2+\frac{25}{3}\mu^3-\frac{13}{3}\lambda^2-\frac{26}{3}\lambda \mu-\frac{13}{3}\mu^2 \geq 0
    \end{align*}
    Once again, labelling the lesser and greater roots of the above quadratic by $r^{-},r^+$, we are able to numerically verify that, under the boundary conditions $0 \leq \lambda,\mu \leq 1/2$, there exists no solution to $\sigma \leq r^{-}$, and there exists a single exact solution for $\sigma \geq r^{+}$, namely $\sigma=-1/2$ which is attained when $\mu=\lambda=1/2$. However, this point attains the value $\Delta(B)=\frac{25}{12}$ exactly, so this does not constitute a contradiction to our claim.
    \\
    We may now assume instead that
    \begin{align}
        1-(1-\lambda-\mu)^2-k(1+\sigma)^2 \leq k(1-\sigma^2) \iff k \geq \frac{1}{2(1+\sigma)}(1-(1-\lambda-\mu)^2). \label{14}
    \end{align}
    Then
    \begin{align*}
        \Delta(B) \leq \left(1+\frac{\lambda^2}{k}\right)\left(1+\frac{\mu^2+k\sigma^2}{k(1-\sigma^2)}\right)=\left(1+\frac{\lambda^2}{k}\right)\left(\frac{1}{1-\sigma^2}+\frac{\mu^2}{k(1-\sigma^2)}\right).
    \end{align*}
    Clearly the upper bound attains its maximum value when $k$ is minimal, so by inequality \ref{14},
    \begin{align*}
        \Delta(B) \leq \left(1+\frac{2\lambda^2(1+\sigma)}{1-(1-\lambda-\mu)^2}\right)\left(\frac{1}{1-\sigma^2}+\frac{2\mu^2}{(1-\sigma)(1-(1-\lambda-\mu)^2)}\right).
    \end{align*}
    However, this inequality is identical to the inequality in \ref{13}, and since we have already shown that this inequality does not rise above $\frac{25}{12}$ in value for $0 \leq \lambda,\mu \leq 1/2$, we are also done for this case.
    \\
    We have confirmed our claim for all $-1/2 \leq \sigma \leq 1/3$, so now we assume that $1/3 \leq \sigma \leq 1/2$. By inequalities \ref{4}, \ref{6}, we have
    \begin{align*}
        l \geq \max\{1-(\lambda-\mu)^2-k(1-\sigma)^2, k(1-\sigma^2)\}.
    \end{align*}
    Assume first that
    \begin{align}
        1-(\lambda-\mu)^2-k(1-\sigma)^2 \geq k(1-\sigma^2) \iff k \leq \frac{1}{2(1-\sigma)}\left(1-(\lambda-\mu)^2\right). \label{15}
    \end{align}
    Then
    \begin{align*}
        \Delta(B) \leq \left(1+\frac{\lambda^2}{k}\right)\left(1+\frac{\mu^2+k\sigma^2}{1-(\lambda-\mu)^2-k(1-\sigma)^2}\right).
    \end{align*}
    Let
    \begin{align*}
        F(k)=\left(1+\frac{\lambda^2}{k}\right)\left(1+\frac{\mu^2+k\sigma^2}{1-(\lambda-\mu)^2-k(1-\sigma)^2}\right).
    \end{align*}
    Then
    \begin{align*}
        F^{\prime \prime}(k)=\frac{\epsilon}{k^3(1-(\lambda-\mu)^2-k(1-\sigma)^2)^3},
    \end{align*}
    where
    \begin{align*}
        \epsilon&=2\lambda^2(\mu^2+k\sigma^2)(1-(\lambda-\mu)^2-k(1-\sigma)^2)^2+2\lambda^2(1-(\lambda-\mu)^2-k(1-\sigma)^2)^3 \\&-2\lambda^2k(1-\sigma)^2(\mu^2+k\sigma^2)(1-(\lambda-\mu)^2-k(1-\sigma)^2)-2\lambda^2k\sigma^2(1-(\lambda-\mu)^2-k(1-\sigma)^2)^2\\&+2\lambda^2k^2(1-\sigma)^4(\mu^2+k\sigma^2)+2k^3(1-\sigma)^4(\mu^2+k\sigma^2)+2\lambda^2k^2\sigma^2(1-\sigma)^2(1-(\lambda-\mu)^2-k(1-\sigma)^2)\\&+2k^3\sigma^2(1-\sigma)^2(1-(\lambda-\mu)^2-k(1-\sigma)^2)
        \\&=\lambda(\mu^2+k\sigma^2)(1-(\lambda-\mu)^2-2k(1-\sigma)^2)^2
        \\&+\lambda^2(1-(\lambda-\mu)^2-k(1-\sigma)^2)(1-(\lambda-\mu)^2-k((1-\sigma)^2+\sigma^2))^2
        \\&+\lambda^2(\mu^2+k\sigma^2)(1-(\lambda-\mu)^2-k(1-\sigma)^2)^2+\lambda^2(1-(\lambda-\mu)^2-k(1-\sigma)^2)^3
        \\&+\lambda^2k^2(1-\sigma)^4(\mu^2+k\sigma^2)+2k^3(1-\sigma)^4(\mu^2+k\sigma^2)
        \\&+\lambda^2k^2\sigma^2(1-\sigma)^2(1-(\lambda-\mu)^2-k(1-\sigma)^2)+k^3(2\sigma^2(1-\sigma)^2-\sigma^4)(1-(\lambda-\mu)^2-k(1-\sigma)^2)
        \\&+2k^3\sigma^2(1-\sigma)^2(1-(\lambda-\mu)^2-k(1-\sigma)^2) \geq 0,
    \end{align*}
    since
    \begin{align*}
        2\sigma^2(1-\sigma)^2-\sigma^4=\sigma^2(2-4\sigma+\sigma^4) > 0,
    \end{align*}
    for all $1/3 \leq \sigma \leq 1/2$, and so $F^{\prime \prime}(k)$ is non-negative for all values of $\lambda,k,\sigma,\mu$ over their respective regions, so the maximum value of $F(k)$ occurs either when $k$ is maximum or when it is minimum. Suppose first that $F(k)$ is maximum when $k$ is minimum. Then
    \begin{align*}
        \Delta(B) \leq \left(1+\frac{\lambda^2}{1-\lambda^2}\right)\left(1+\frac{\mu^2+(1-\lambda^2)\sigma^2}{1-(\lambda-\mu)^2-(1-\lambda^2)(1-\sigma)^2}\right).
    \end{align*}
    Assume that
    \begin{align*}
        \left(1+\frac{\lambda^2}{1-\lambda^2}\right)\left(1+\frac{\mu^2+(1-\lambda^2)\sigma^2}{1-(\lambda-\mu)^2-(1-\lambda^2)(1-\sigma)^2}\right) \geq \frac{25}{12},
    \end{align*}
    which gives us the quadratic
    \begin{align*}
        &\frac{25}{12}(1-\lambda^2)^2\sigma^2+\left(2(1-\lambda^2)-\frac{25}{6}(1-\lambda^2)^2\right)\sigma\\&+1-(\lambda-\mu)^2-\frac{37}{12}(1-\lambda^2)+\mu^2+\frac{25}{12}(1-\lambda^2)(\lambda-\mu)^2+\frac{25}{12}(1-\lambda^2)^2 \geq 0.
    \end{align*}
    Let $r^{-},r^+$ respectively denote the lesser and greater roots of the quadratic above. By inequality \ref{15}, we obtain the boundary conditions
    \begin{align}
        &0 \leq \lambda \leq 1/2, \label{16}
        \\&0 \leq \mu \leq 2\lambda. \label{17}
    \end{align}
    Again, we need either $\sigma \leq r^{-}$ or $\sigma \geq r^+$. However, we can numerically verify that under boundary conditions \ref{16},\ref{17}, the solution $\sigma \leq r^{-}$ requires $\sigma <1/3$ and the solution $\sigma \geq r^+$ requires that $\sigma >1/2$, both of which are clearly contradictions under the assumptions made.
    \\
    Now, assume that $F(k)$ attains its maximum when $k$ is maximum. Then
    \begin{align}
        \Delta(B) \leq \left(1+\frac{2\lambda^2(1-\sigma)}{1-(\lambda-\mu)^2}\right)\left(1+\frac{\mu^2+\frac{\sigma^2}{2(1-\sigma)}(1-(\lambda-\mu)^2}{\frac{1}{2}(1+\sigma)(1-(\lambda-\mu)^2)}\right). \label{18}
    \end{align}
    Assume that
    \begin{align*}
        \left(1+\frac{2\lambda^2(1-\sigma)}{1-(\lambda-\mu)^2}\right)\left(1+\frac{\mu^2+\frac{\sigma^2}{2(1-\sigma)}(1-(\lambda-\mu)^2}{\frac{1}{2}(1+\sigma)(1-(\lambda-\mu)^2)}\right) \geq \frac{25}{12}.
    \end{align*}
    Rearranging gives us the quadratic
    \begin{align*}
        &\frac{25\lambda^4-100\lambda^3\mu+198\lambda^2\mu^2-100\lambda\mu^3+25\mu^4-50\lambda^2+100\lambda\mu-50\mu^2+25}{12}\sigma^2
        \\&+2(\lambda^4-2\lambda^3\mu-2\lambda^2\mu^2-2\lambda\mu^3+\mu^4-\lambda^2-\mu^2)\sigma
        \\&-\frac{37}{12}\lambda^4+\frac{25}{3}\lambda^3\mu-\frac{13}{2}\lambda^2\mu^2+\frac{25}{3}\lambda\mu^3-\frac{37}{12}\mu^4+\frac{25}{6}\lambda^2-\frac{13}{3}\lambda \mu+\frac{25}{6}\mu^2-\frac{13}{12} \geq 0.
    \end{align*}
    Again, letting $r^{-},r^+$ denote the lesser and greater roots of the quadratic above, we verify numerically that $r^{-} < 1/3$ and $r^+ \geq 1/2$ for all $0 \leq \mu, \lambda \leq 1/2$, with equality $\sigma=1/2$ only when $\lambda=\mu=1/2$. We have $\Delta(B)=\frac{25}{12}$ at the point $\sigma=\mu=\lambda=1/2$ but does not exceed $\frac{25}{12}$ at any other point, so this does not constitute a contradiction to our claim.
    \\
    Finally, assume that
    \begin{align*}
        1-(\lambda-\mu)^2-k(1-\sigma)^2 \leq k(1-\sigma^2) \iff k \geq \frac{1}{2(1-\sigma)}\left(1-(\lambda-\mu)^2\right),
    \end{align*}
    so
    \begin{align*}
        \Delta(B) \leq \left(1+\frac{\lambda^2}{k}\right)\left(1+\frac{\mu^2+k\sigma^2}{k(1-\sigma^2)}\right)=\left(\frac{1}{1-\sigma^2}+\frac{\mu^2}{k(1-\sigma^2)}\right).
    \end{align*}
    Clearly the upper bound attains its maximum when $k$ is minimal, so
    \begin{align*}
        \Delta(B) \leq \left(1+\frac{2\lambda^2(1-\sigma)}{1-(\lambda-\mu)^2}\right)\left(\frac{1}{1-\sigma^2}+\frac{2\mu^2}{(1+\sigma)(1-(\lambda-\mu)^2)}\right).
    \end{align*}
    However, the function determining the upper bound for $\Delta(B)$ here is identical to that in \ref{18}, for which we have verified is less than or equal to $\frac{25}{12}$ over all values of $\sigma, \mu, \lambda$ in their relative intervals, so we do not need to deal with this case. Since this exhausts all possible values for all the variables, we have verified that $\Delta_3=\frac{25}{12}$, and the quadratic form
    \begin{align*}
        \left(x+\frac{1}{2}(y+z)\right)^2+\left(y \pm \frac{1}{2}z\right)^2+\frac{3}{4}z^2
    \end{align*}
    attains this upper bound, so this is the best possible upper bound on the orthogonality defect for three dimensional HKZ reduced bases.
\section{Proof of Theorem \ref{thm4}}
    Suppose that $B=\{\mathbf{b}_1,\mathbf{b}_2,\dots,\mathbf{b}_n\}$ is the basis for a rank $n \geq 4$ lattice, and assume that $B$ is HKZ reduced. Then note that the subbasis $\{\mathbf{b}_1,\mathbf{b}_2,\mathbf{b}_3\}$ must also be HKZ reduced, and so
    \begin{align*}
        \frac{\|\mathbf{b}_1\|^2\|\mathbf{b}_2\|^2\|\mathbf{b}_3\|^2}{\|\mathbf{b}_1\|^2\|\mathbf{b}_2(2)\|^2\|\mathbf{b}_3(3)\|^2} \leq \frac{25}{12}.
    \end{align*}
    Then
    \begin{align*}
        \Delta(B)=\prod_{i=1}^n \frac{\|\mathbf{b}_i\|^2}{\|\mathbf{b}_i(i)\|^2} \leq \frac{25}{12}\prod_{i=4}^n \frac{\|\mathbf{b}_i\|^2}{\|\mathbf{b}_i(i)\|^2}.
    \end{align*}
    Also note that the lattice generated by the basis $\{\mathbf{b}_4(4),\mathbf{b}_5(4),\dots,\mathbf{b}_n(4)\}$ is HKZ reduced, by definition. For any $i \geq 4$, we have
    \begin{align*}
        \|\mathbf{b}_i\|^2=\|\mathbf{b}_i(4)\|^2+\sum_{j=1}^3|\mu_{i,j}|^2\|\mathbf{b}_j(j)\|^2 \leq \|\mathbf{b}_i(4)\|^2+\frac{1}{4}\sum_{j=1}^3\|\mathbf{b}_j(j)\|^2.
    \end{align*}
    Using Proposition \ref{propn1},
    \begin{align*}
        \|\mathbf{b}_1\|^2 \leq 2\|\mathbf{b}_4(4)\|^2, \|\mathbf{b}_2(2)\|^2 \leq \frac{3}{2}\|\mathbf{b}_4(4)\|^2, \|\mathbf{b}_3(3)\|^2 \leq \frac{4}{3}\|\mathbf{b}_4(4)\|^2,
    \end{align*}
    so
    \begin{align*}
    \|\mathbf{b}_i\|^2 \leq \|\mathbf{b}_i(4)\|^2 +\frac{1}{4}\left(2+\frac{3}{2}+\frac{4}{3}\right)\|\mathbf{b}_4(4)\|^2=\|\mathbf{b}_i(4)\|^2+\frac{29}{24}\|\mathbf{b}_4(4)\|^2.
    \end{align*}
    Lagarias et. al. proved that if $B$ is HKZ reduced, then
    \begin{align*}
        \|\mathbf{b}_i(i)\|^2 \leq \lambda_i^2,
    \end{align*}
    for all $1 \leq i \leq n$, and so for all $i \geq 4$, letting $\lambda_i(4)$ denote the $i$th successive minima in the lattice generated by $\{\mathbf{b}_4(4),\mathbf{b}_5(4),\dots,\mathbf{b}_n(4)\}$,
    \begin{align*}
        \|\mathbf{b}_i\|^2 &\leq \|\mathbf{b}_i(4)\|^2+\frac{29}{24}\|\mathbf{b}_4(4)\|^2 =\|\mathbf{b}_i(i)\|^2+\sum_{j=4}^i |\mu_{i,j}|^2\|\mathbf{b}_j(j)\|^2+\frac{29}{24}\|\mathbf{b}_4(4)\|^2 \\&\leq \lambda_i(4)^2+\frac{1}{4}\sum_{j=4}^i \lambda_j(4)^2+\frac{29}{24}\lambda_4(4)^2 \leq \left(\frac{i}{4}+\frac{29}{24}\right)\lambda_i(4)^2.
    \end{align*}
    Then, using Minkowski's second theorem we get
    \begin{align*}
        \Delta(B) \leq \frac{25}{12}\prod_{i=4}^n\frac{\|\mathbf{b}_i\|^2}{\|\mathbf{b}_i(i)\|^2} \leq \frac{25}{12}\prod_{i=4}^n\left(\frac{i}{4}+\frac{29}{24}\right)\frac{\lambda_i(4)^2}{\|\mathbf{b}_i(i)\|^2} \leq \frac{25}{12}\gamma_{n-3}^{n-3} \prod_{i=4}^n \left(\frac{i}{4}+\frac{29}{24}\right),
    \end{align*}
    as required.
    \section{Summary and Future Work}
    In this work, we determined the exact maximal value that the orthogonality defect can take for HKZ reduced bases up to dimension $3$, and used this result to calculate an upper bound on the orthogonality defect of HKZ bases of arbitrary rank. The latter result seems to be sharper than any other bounds on the orthogonality defect for HKZ reduced bases in existing literature, according to the authors' knowledge.
    \\
    In the future, we will investigate whether the techniques applied in this paper can be used to sharpen the bounds on the orthogonality defect further, and to determine the value of $\Delta_n$ exactly for $n \geq 4$. We conjecture that $\Delta_n=\gamma_n^n$ for all $n \geq 4$, and preliminary computational results seem to agree with this conjecture.
    
\end{document}